\documentclass[a4paper,12pt]{article}

\usepackage[margin=25mm]{geometry}

\usepackage[latin1]{inputenc}
\usepackage[T1]{fontenc}
\usepackage[english]{babel}

\usepackage{amsmath,amssymb}
\usepackage{bbm}

\newcommand*{\bs}{\boldsymbol}

\newcommand*{\mrm}{\mathrm}

\newcommand*{\diff}{\mathop{}\!\mathrm{d}}
\newcommand*{\pd}{\partial}
\newcommand*{\im}{\ensuremath{\mathrm{i}}}
\newcommand*{\e}{\mathop{\mathrm{e}}\nolimits}
\newcommand*{\CC}{\mathbb{C}}
\newcommand*{\Z}{\mathbb{Z}}
\newcommand*{\defin}{\stackrel{\mathrm{def}}{=}}
\DeclareMathOperator{\Real}{Re}
\DeclareMathOperator{\Imag}{Im}
\DeclareMathOperator*{\Res}{Res}
\newcommand*{\diag}{\mathrm{diag}}
\DeclareMathOperator{\sign}{sign}

\newcommand*{\speed}{\mathrm{s}}
\newcommand*{\vol}{\mathrm{v}}

\def\Xint#1{\mathchoice
  {\XXint\displaystyle\textstyle{#1}}%
  {\XXint\textstyle\scriptstyle{#1}}%
  {\XXint\scriptstyle\scriptscriptstyle{#1}}%
  {\XXint\scriptscriptstyle\scriptscriptstyle{#1}}%
  \!\int}
\def\XXint#1#2#3{{\setbox0=\hbox{$#1{#2#3}{\int}$}
    \vcenter{\hbox{$#2#3$}}\kern-.5\wd0}}
\def\pvint{\Xint-}

\numberwithin{equation}{section}
\allowdisplaybreaks[4]

\usepackage[numbers, sort&compress]{natbib}

\providecommand{\url}[1]{\texttt{#1}}

\usepackage[nottoc]{tocbibind}

\newcommand*{\ie}{i.\,e.}
\newcommand*{\eg}{e.\,g.}

\usepackage{indentfirst}
\usepackage{graphicx}
\usepackage[font=small,labelfont=bf]{caption}
\usepackage{xcolor}

\usepackage[unicode,colorlinks]{hyperref}

\definecolor{darkred}{rgb}{0.7,0,0}
\definecolor{darkgreen}{rgb}{0,0.5,0}
\definecolor{darkblue}{rgb}{0,0,0.8}

\hypersetup{%
  bookmarksopen=true,
  bookmarksopenlevel=2,
  bookmarksnumbered=true,
  linkcolor=darkblue,
  citecolor=darkblue,
  urlcolor=darkblue,
  pdfpagelayout=OneColumn,
  pdfstartview={XYZ null null 1},
  pdfview={XYZ null null null}
}

\begin{document}

\title{%
  \bfseries
  Integral representation of solutions to initial-boundary value problems in the framework of the Guyer--Krumhansl heat equation
}

\author{%
  \bfseries
  \large
  Sergey~A.~Rukolaine\thanks{E-mail address: \texttt{rukol@ammp.ioffe.ru}}
}

\date{%
  \small
  Ioffe Institute, 26 Polytekhnicheskaya, St.\,Petersburg 194021, Russia
}

\maketitle

\begin{abstract}
  We consider initial-boundary value problems (IBVPs) on a finite interval for the system of the energy balance equation and Guyer--Krumhansl constitutive equation. Boundary conditions comprise various models of behavior of a physical system at the boundaries, including boundary conditions describing Newton's law, which states that the heat flux at the boundary is directly proportional to the difference in the temperature of the physical system and ambient temperature. In this case the boundary conditions express the relationship of unknown functions (temperature or internal energy and heat flux) with each other. To solve the problems, we apply the Fokas unified transform method. To illustrate the final formulae, we consider a numerical example for a special case in which heat exchange at one of the boundaries obeys Newton's law.
\end{abstract}

\section{Introduction}

The Guyer--Krumhansl heat equation is obtained from the energy balance equation and Guyer--Krumhansl constitutive equation. The energy balance equation is expressed by the equation
\begin{equation}
  \label{eq:EnergyEqSource}
  C_\vol^{} \pd_t^{} T
  + \nabla \cdot \bs{q}
  =
  f(\bs{r},t),
\end{equation}
where $T \equiv T(\bs{r},t)$ is temperature, $C_\vol^{}$ is volumetric heat capacity, $\bs{q} \equiv \bs{q}(\bs{r},t)$ is the vector of heat flux, $f$ is the heat source term.

The Guyer--Krumhansl constitutive equation, relating  heat flux and temperature, is given by \cite{GuyerKrumhansl:1966a, JosephPreziosi:1989, JouEtAl:2010}
\begin{equation}
  \label{eq:GuKrRelOrig}
  \pd_t^{} \bs{q}
  + \frac{1}{\tau_{\mrm{R}}^{}} \bs{q}
  + \frac{c_\speed^2 C_\vol^{}}{3} \nabla T
  =
  \frac{\tau_{\mrm{N}}^{} c_\speed^2}{5} \left( \Delta \bs{q} + 2 \nabla \nabla \cdot \bs{q} \right),
\end{equation}
where $c_\speed^{}$ is the Debye (average) speed of phonons, $\tau_{\mrm{N}}^{}$ is a
relaxation time  for momentum-conserving collisions (normal (N-)processes)
and $\tau_{\mrm{R}}^{}$ is a relaxation time for momentum-nonconserving
collisions (resistive (R-)processes) in the phonon gas.
The Guyer--Krum\-hansl constitutive equation~\eqref{eq:GuKrRelOrig} can be written as
\begin{equation}
  \label{eq:GKRel}
  \tau \pd_t^{} \bs{q}
  + \bs{q}
  + \kappa \nabla T
  =
  \ell^2 \left( \Delta \bs{q} + 2 \nabla \nabla \cdot \bs{q} \right),
\end{equation}
where
\begin{equation*}
  \tau
  =
  \tau_{\mrm{R}}^{},
  \quad
  \kappa
  =
  \frac{\tau_{\mrm{R}}^{} c_\speed^2 C_\vol^{}}{3}
  \equiv
  \frac{\ell^2 C_\vol^{}}{3 \tau_{\mrm{R}}^{}},
  \quad
  \ell^2
  =
  \frac{\tau_{\mrm{R}}^{} \tau_{\mrm{N}}^{} c_\speed^2}{5},
\end{equation*}
$\tau$ is the relaxation time, $\kappa$ is the thermal conductivity, $\ell = c_\speed^{} \tau$ is the mean free path of the heat carriers (phonons).

It follows from the energy balance equation~\eqref{eq:EnergyEqSource} and the Guyer--Krumhansl constitutive equation~\eqref{eq:GKRel} that temperature satisfies the equation
\begin{equation}
  \label{eq:GKEqGen}
  C_\vol^{} \left[
    \tau \pd_t^2 T
    + \pd_t^{} T
    - \mu^2 \pd_t^{} \Delta T
    - \alpha \Delta T
  \right]
  =
  \tau \pd_t^{} f
  + f
  - \mu^2 \Delta f.
\end{equation}
where
\begin{equation*}
  \alpha
  =
  \frac{\kappa}{C_\vol^{}},
  \quad
  \mu^2
  =
  3 \ell^2,
\end{equation*}
$\alpha$ is thermal diffusivity.

In the absence of the heat source, \ie, if $f=0$, Eq.~\eqref{eq:GKEqGen} takes the form
\begin{equation}
  \label{eq:GKEqNoSource}
  \tau \pd_t^2 T
  + \pd_t^{} T
  - \mu^2 \pd_t^{} \Delta T
  - \alpha \Delta T
  = 0,
\end{equation}
This is the Guyer--Krumhansl heat equation, as it is usually written~\cite{JosephPreziosi:1989}.

The Guyer--Krumhansl model was applied to heat conduction in Refs.~\cite{LebonEtAl:2011PRSA, VanEtAl:2017EPL, Kovacs:2018IJHMT, RiethEtAl:2018IJHMT, FeherKovacs:2021, Xu:2021PRSA, RamosEtAl:2023AMM, Toth:2023IJHMT, TothEtAl:2024CMT, Toth:2025IJHMT, Tur-PratsEtAl:2024IJHMT}.

Solving initial value problems (in the entire space) for this and other equations obtained in the framework of extended thermodynamics and analysing their solutions is straightforward \cite{RukolaineSamsonov:2013, Rukolaine:2024ZAMP, Rukolaine:2024TP}.
In the papers~\cite{Kovacs:2018IJHMT, RiethEtAl:2018IJHMT, FeherKovacs:2021, RamosEtAl:2023AMM, Toth:2023IJHMT, Toth:2025IJHMT} initial-boundary value problems (IBVPs) for the Guyer--Krumhansl heat equation were solved both analytically and numerically for particular boundary conditions.
However, solving IBVPs analytically in general case may encounter difficulties.

Below, for simplicity, we will consider equations in one-dimensional space.

In one-dimensional space Eq.~\eqref{eq:GKEqNoSource} takes the form
\begin{equation}
  \label{eq:GKEqOneDimNoSource}
  \tau \pd_t^2 T
  + \pd_t^{} T
  - \mu^2 \pd_t^{} \pd_x^2 T
  - \alpha \pd_x^2 T
  = 0.
\end{equation}
If both boundary conditions set the temperature behavior, then the IBVP can easily be solved by the Fourier method or Laplace transform. If both boundary conditions set the heat flux behavior, then the analytical solution for temperature can also be easily found. Indeed, heat flux satisfies the same equation as temperature:
\begin{equation}
  \label{eq:GKEqFluxOneDim}
  \tau \pd_t^2 q
  + \pd_t^{} q
  - \mu^2 \pd_t^{} \pd_x^2 q
  - \alpha \pd_x^2 q
  = 0.
\end{equation}
First, the IBVP for this equation is solved, and then temperature is determined using the energy balance equation~\eqref{eq:EnergyEqSource} (in one dimension).

However, other boundary conditions may be set as well, \eg,
one or both boundary conditions can describe Newton's law (convection boundary conditions), which states that the heat flux at the boundary is directly proportional to the difference in the temperature of the physical system and ambient temperature. In this case IBVPs cannot be formulated for either equation~\eqref{eq:GKEqOneDimNoSource} or equation~\eqref{eq:GKEqFluxOneDim}. In this case IBVPs should be formulated for the system of the energy balance equation~\eqref{eq:EnergyEqSource} and Guyer--Krumhansl constitutive equation~\eqref{eq:GKRel} (both in one dimension if we stay in one-dimensional space). And in this case, the Fourier method is not applicable. Application of the Laplace transform to these problems is problematic as well.
All these problems can be solved by the Fokas unified transform method (UTM) \cite{Fokas:1997PRSA, Fokas:2002IMAJAM, FokasPelloni:2005IMAJAM, FokasPelloni:2005MPAG, Fokas:2008, MantzavinosFokas:2013, DeconinckEtAl:2014SIAMRev, DeconinckEtAl:2018QAM, JohnstonEtAl:2021SAM, FokasKalimeris:2022IMAJAM, MantzavinosMitsotakis:2023JMPA, ChatziafratisEtAl:2024ZAMP, ChatziafratisEtAl:2024MMAS, ChatziafratisEtAl:2024ZAMM}.
In this paper we solve these IBVPs by the UTM.
The UTM was applied to systems of equations in Refs.~\cite{FokasPelloni:2005MPAG, DeconinckEtAl:2018QAM, JohnstonEtAl:2021SAM, MantzavinosMitsotakis:2023JMPA, ChatziafratisEtAl:2024MMAS, ChatziafratisEtAl:2024ZAMP}.
The peculiarity of the problems considered in this paper is that the boundary conditions generally express the relationship of unknown functions (temperature and heat flux) with each other.

\section{Statement of the initial-boundary value problem}
\label{sec:IBVP}

We consider the IBVP for the system of the energy balance equation (including the heat source term) and Guyer--Krumhansl constitutive equation on the interval $(0,l)$. For convenience, we consider internal energy $e = C_\vol^{} T$ instead of temperature $T$. Therefore, we solve the system of the equations
\begin{subequations}
  \label{eq:IBVPSystem}
  \begin{align}
    \label{eq:IBVPEqA}
    &\pd_t^{} e + \pd_x^{} q
      =
      f(x,t),
      \quad
      x \in (0,l),
      \quad
      t > 0,
    \\[0.5ex]
    \label{eq:IBVPEqB}
    &\tau \pd_t^{} q
      + q
      - \mu^2 \pd_x^2 q
      + \alpha \pd_x^{} e
      = 0,
      \quad
      x \in (0,l),
      \quad
      t > 0,
  \end{align}
\end{subequations}
where $f$ is the heat source term.
Initial conditions for the system~\eqref{eq:IBVPSystem} set the initial state of the physical system and are given by
\begin{equation}
  \label{eq:IBVPInitCond}
  \left. e \right|_{t=0}^{}
  =
  \varphi(x),
  \quad
  \left. q \right|_{t=0}^{}
  =
  \psi(x).
\end{equation}
The boundary conditions are given by
\begin{equation}
  \label{eq:IBVPBoundCond}
  \left. (\gamma_0^{} e + \delta_0^{} q) \right|_{x=0}^{}
  =
  g(t),
  \quad
  \left. (\gamma_l^{} e + \delta_l^{} q) \right|_{x=l}^{}
  =
  h(t),
\end{equation}
where
$|\gamma_0^{}| + |\delta_0^{}| > 0$ and $|\gamma_l^{}| + |\delta_l^{}| > 0$.

It is convenient to write the system~\eqref{eq:IBVPSystem} in the matrix form
\begin{equation}
  \label{eq:IBVPSystemMatrixForm}
  \pd_t^{} \bs{u}
  + \varLambda(-\im \pd_x^{}) \bs{u}
  = \bs{F}(x,t),
  \quad
  \bs{u}
  =
  \begin{pmatrix}
    e \\ q
  \end{pmatrix},
  \quad
  \bs{F}(x,t)
  =
  \begin{pmatrix}
    f(x,t) \\ 0
  \end{pmatrix},
\end{equation}
where the matrix $\varLambda$ is given by
\begin{equation*}
  \varLambda(k)
  =
  \begin{pmatrix}
    0 & \im k
    \\[5pt]
    \dfrac{\alpha}{\tau} \im k & \dfrac{\mu^2}{\tau} k^2 + \dfrac{1}{\tau}
  \end{pmatrix}
  \equiv
  \begin{pmatrix}
    0 & \im k
    \\[5pt]
    \im k \beta & \theta k^2 + \dfrac{1}{\tau}
  \end{pmatrix}.
\end{equation*}
The initial conditions~\eqref{eq:IBVPInitCond} take the form
\begin{equation*}
  \left. \bs{u} \right|_{t=0}^{}
  =
  \bs{u}_0(x)
  \equiv
  \begin{pmatrix}
    \varphi(x) \\[2pt] \psi(x)
  \end{pmatrix}.
\end{equation*}

\section{Solving the problem}
\label{sec:Solving}

\subsection{The local relations}

First we need to find dispersion relations. They are determined by the equation
\begin{equation*}
  \det \left[ \varLambda(k) - \omega I \right]
  = 0,
\end{equation*}
or, equivalently,
\begin{equation*}
  \omega^2
  - \frac{1 + \mu^2 k^2}{\tau} \omega
  + \frac{\alpha}{\tau} k^2
  = 0.
\end{equation*}
The roots of the equation, given by
\begin{equation}
  \label{eq:omegaOneTwo}
  \omega_{1,2}(k)
  =
  \frac{1 + \mu^2 k^2 \pm \sqrt{(1 + \mu^2 k^2)^2 - 4 \alpha \tau k^2}}{2\tau},
\end{equation}
determine the dispersion relations.

Next, we need to write the system~\eqref{eq:IBVPSystemMatrixForm} in divergence form, which is called the local relation. We have the chain of equalities
\begin{multline*}
  \pd_t^{} \!\left( \e^{-\im k x} \e^{\varLambda(k) t} \bs{u} \right)
  =
  \e^{-\im k x} \e^{\varLambda(k) t} \left[
    \varLambda(k) \bs{u}
    + \pd_t^{} \bs{u}
  \right]
  \\
  =
  \e^{-\im k x} \e^{\varLambda(k) t} \left[
    \varLambda(k) \bs{u}
    - \varLambda(-\im \pd_x^{}) \bs{u}
    + \bs{F}(x,t)
  \right]
  \\
  =
  \e^{-\im k x} \e^{\varLambda(k) t} \left[
    (k + \im \pd_x^{}) \left. \frac{\varLambda(k) - \varLambda(l)}{k - l} \right|_{l = -\im \pd_x^{}} \bs{u}
    + \bs{F}(x,t)
  \right]
  \\
  =
  \pd_x^{} \!\left( \e^{-\im k x} \e^{\varLambda(k) t} X(k) \bs{u} \right)
  + \e^{-\im k x} \e^{\varLambda(k) t} \bs{F}(x,t),
\end{multline*}
where $X$ is a matrix operator, given by
\begin{equation*}
  X(k)
  =
  \im \left. \frac{\varLambda(k) - \varLambda(l)}{k - l} \right|_{l = -\im \pd_x^{}}^{}.
\end{equation*}
So we get the local relation
\begin{equation}
  \label{eq:LocRel}
  \pd_t^{} \!\left( \e^{-\im k x} \e^{\varLambda(k) t} \bs{u} \right)
  - \pd_x^{} \!\left( \e^{-\im k x} \e^{\varLambda(k) t} X(k) \bs{u} \right)
  =
  \e^{-\im k x} \e^{\varLambda(k) t} \bs{F}(x,t).
\end{equation}
In this problem
\begin{equation*}
  X(k)
  =
  \begin{pmatrix}
    0 & -1
    \\[3pt]
    - \dfrac{\alpha}{\tau} & \dfrac{\mu^2}{\tau} (\im k + \pd_x^{})
  \end{pmatrix}
  \equiv
  \begin{pmatrix}
    0 & -1
    \\[2pt]
    - \beta & \theta (\im k + \pd_x^{})
  \end{pmatrix},
  \quad
  \beta = \dfrac{\alpha}{\tau},
  \quad
  \theta = \dfrac{\mu^2}{\tau}.
\end{equation*}

It is convenient to diagonalize the matrix $\varLambda$ (see Refs.~\cite{FokasPelloni:2005MPAG, DeconinckEtAl:2018QAM}):
\begin{equation*}
  \varLambda(k)
  =
  S(k) \varOmega(k) S^{-1}(k),
\end{equation*}
where $\varOmega$ is a diagonal matrix, given by
\begin{equation*}
  \varOmega(k)
  =
  \diag \big( \omega_1^{}(k), \omega_2^{}(k) \big),
\end{equation*}
so that
\begin{equation*}
  \e^{-\varLambda(k) t}
  =
  S(k) \e^{-\varOmega(k) t} S^{-1}(k).
\end{equation*}
In this case the local relation~\eqref{eq:LocRel} can be written in the form
\begin{equation}
  \label{eq:LocRelM}
  \pd_t^{} \!\left( \e^{-\im k x} \e^{\varOmega(k) t} S^{-1}(k) \bs{u} \right)
  - \pd_x^{} \!\left( \e^{-\im k x} \e^{\varOmega(k) t} S^{-1}(k) X(k) \bs{u} \right)
  =
  \e^{-\im k x} \e^{\varOmega(k) t} S^{-1}(k) \bs{F}(x,t).
\end{equation}

The columns of the matrix $S$ are the eigenvectors of the matrix $\varLambda$, and in this problem
\begin{equation*}
  S(k)
  =
  \frac{1}{\omega_1^{} - \omega_2^{}}
  \begin{pmatrix}
    1 & 1
    \\
    \omega_1^{} / \im k & \omega_2^{} / \im k
  \end{pmatrix},
  \qquad
  S^{-1}(k)
  =
  \begin{pmatrix}
    -\omega_2^{} & \im k
    \\
    \omega_1^{} & -\im k
  \end{pmatrix}.
\end{equation*}

\subsection{The global relations}

Integrating the local relations~\eqref{eq:LocRel} and \eqref{eq:LocRelM} over the domain $\{(x,t) \in (0,l) \times (0,T) \}$, $T > 0$, using Green's theorem, and replacing $T$ with $t$, we obtain the global relations
\begin{multline}
  \label{eq:GlobRel}
  \hat{\bs{u}}(k,t)
  - \e^{-\Lambda(k) t} \hat{\bs{u}}_0^{}(k)
  + S(k) \bs{a}(k,t)
  - \e^{-\im k l} S(k) \bs{b}(k,t)
  =
  \int_0^t \e^{-\Lambda(k) (t-s)} \hat{\bs{F}}(k,s) \diff s,
  \\
  \quad
  k \in \CC,
\end{multline}
and
\begin{equation}
  \label{eq:GlobRelM}
  S^{-1}(k) \hat{\bs{u}}(k,t)
  - \e^{-\varOmega(k) t} S^{-1}(k) \hat{\bs{u}}_0^{}(k)
  + \bs{a}(k,t)
  - \e^{-\im k l} \bs{b}(k,t)
  =
  \tilde{\bs{F}}(k,t),
  \quad
  k \in \CC,
\end{equation}
respectively, where $\hat{\bs{u}}$ and $\hat{\bs{u}}_0^{}$ are the Fourier transforms of the solution $\bs{u}$ and its initial state $\bs{u}_0^{}$, respectively:
\begin{equation}
  \label{eq:HatU}
  \hat{\bs{u}}(k,t)
  =
  \begin{pmatrix}
    \hat{e}(k,t) \\ \hat{q}(k,t)
  \end{pmatrix},
  \quad
  \hat{e}(k,t)
  =
  \int_0^l \e^{-\im k x} e(x,t) \diff x,
  \quad
  \hat{q}(k,t)
  =
  \int_0^l \e^{-\im k x} q(x,t) \diff x,
\end{equation}
and
\begin{equation}
  \label{eq:HatUZero}
  \hat{\bs{u}}_0^{}(k)
  =
  \begin{pmatrix}
    \hat\varphi(k) \\ \hat\psi(k)
  \end{pmatrix},
  \quad
  \hat\varphi(k)
  =
  \int_0^l \e^{-\im k x} \varphi(x) \diff x,
  \quad
  \hat\psi(k)
  =
  \int_0^l \e^{-\im k x} \psi(x) \diff x,
\end{equation}
(we assume that all the functions outside the interval $[0,l]$ are zero, therefore, $\hat{}$ means the Fourier transform),
the vectors $\bs{a}$ and $\bs{b}$ are expressed through the boundary values of the solution $\bs{u}$ by
\begin{subequations}
  \label{eq:ab}
  \begin{align}
    \label{eq:a}
    \bs{a}(k,t)
    &=
      \int_0^t \e^{-\varOmega(k) (t-s)} S^{-1}(k) X \bs{u}(0,s) \diff s
    \\
    \notag
    &=
    \begin{pmatrix}
      - \big[
      \im k \beta \tilde{e}_0^{(1)}
      - (\omega_2^{} - \theta k^2) \tilde{q}_0^{(1)}
      - \im k \theta \tilde{q}_1^{(1)}
      \big] (0,t)
      \\[3pt]
      \big[
      \im k \beta \tilde{e}_0^{(2)}
      - (\omega_1^{} - \theta k^2) \tilde{q}_0^{(2)}
      - \im k \theta \tilde{q}_1^{(2)}
      \big] (0,t)
    \end{pmatrix},
    \\[2pt]
    \label{eq:b}
    \bs{b}(k,t)
    &=
    \int_0^t \e^{-\varOmega(k) (t-s)} S^{-1}(k) X \bs{u}(l,s) \diff s
    \\
    \notag
    &=
    \begin{pmatrix}
      - \big[
      \im k \beta \tilde{e}_0^{(1)}
      - (\omega_2^{} - \theta k^2) \tilde{q}_0^{(1)}
      - \im k \theta \tilde{q}_1^{(1)}
      \big] (l,t)
      \\[3pt]
      \big[
      \im k \beta \tilde{e}_0^{(2)}
      - (\omega_1^{} - \theta k^2) \tilde{q}_0^{(2)}
      - \im k \theta \tilde{q}_1^{(2)}
      \big] (l,t)
    \end{pmatrix},
  \end{align}
\end{subequations}
the functions $\tilde{e}_j^{(1,2)}$ and $\tilde{q}_j^{(1,2)}$ are time transforms of $e$ and $q$ and the derivatives $\pd_x^{} e$ and $\pd_x^{} q$, respectively, defined by
\begin{subequations}
  \label{eq:tildeeq}
  \begin{align}
    \label{eq:tildee}
    \tilde{e}_j^{(m)}(x,t)
    \equiv
    \tilde{e}_j^{}(\omega_m^{}(k),x,t)
    =
    \int_0^t \e^{-\omega_m^{}(k) (t-s)} \pd_x^j e(x,s) \diff s,
    \\[2pt]
    \label{eq:tildeq}
    \tilde{q}_j^{(m)}(x,t)
    \equiv
    \tilde{q}_j^{}(\omega_m^{}(k),x,t)
    =
    \int_0^t \e^{-\omega_m^{}(k) (t-s)} \pd_x^j q(x,s) \diff s,
  \end{align}
\end{subequations}
(note that $\tilde{e}_j^{(m)}$ and $\tilde{q}_j^{(m)}$ do not change under the transformation $k \to -k$),
$\tilde{\bs{F}}$ is expressed through the Fourier transform of $f$ by
\begin{multline}
  \label{eq:tildeF}
  \tilde{\bs{F}}(k,t)
  \equiv
  \begin{pmatrix}
    \tilde{F}_1^{}(k,t) \\[2pt] \tilde{F}_2^{}(k,t)
  \end{pmatrix}
  =
  \int_0^t \e^{-\varOmega(k) (t-s)} S^{-1}(k) \hat{\bs{F}}(k,s) \diff s
  \\
  =
  \int_0^t
  \begin{pmatrix}
    -\omega_2^{} \e^{-\omega_1^{} (t-s)}
    \\[2pt]
    \omega_1^{} \e^{-\omega_2^{} (t-s)}
  \end{pmatrix}
  \hat{f}(k,s) \diff s,
\end{multline}
where
\begin{equation}
  \label{eq:HatF}
  \hat{\bs{F}}(k,t)
  =
  \begin{pmatrix}
    \hat{f}(k,t) \\ 0
  \end{pmatrix},
  \quad
  \hat{f}(k,t)
  =
  \int_0^l \e^{-\im k x} f(x,t) \diff x,
\end{equation}
$\hat{f}$ is the Fourier transform of $f$.

It is important to note that the global relations~\eqref{eq:GlobRel} and \eqref{eq:GlobRelM} are valid for $k \in \CC$. Indeed, the spatial variable $x$ belongs to the finite interval $[0,l]$, therefore, if we assume sufficient decay of the functions $e$, $q$, $\varphi$, $\psi$ and $f$ as $x \to \infty$ for all $t$, then the Fourier transforms can be analytically continued from the real axis of the spectral parameter $k$ onto the entire complex $k$-plane due to the exponential decay there.

\subsection{A preliminary expression for the solution}

We assume that $e$, $q$, $\varphi$, $\psi$ and $f$ are zero, when $x < 0$ or $x > l$. Therefore, the transform $\hat{}$ in Eqs.~\eqref{eq:HatU}, \eqref{eq:HatUZero} and \eqref{eq:HatF} is the Fourier transform over the entire real axis. Applying the inverse Fourier transform to the global relation~\eqref{eq:GlobRel}, we obtain the preliminary integral expression for the solution
\begin{multline}
  \label{eq:SolutionRaw}
  \bs{u}(x,t)
  =
  \frac{1}{2\pi} \int_{-\infty}^\infty \e^{\im k x} \e^{-\varLambda(k) t} \hat{\bs{u}}_0^{}(k) \diff k
  + \frac{1}{2\pi} \int_{-\infty}^\infty \e^{\im k x} \int_0^t \e^{-\varLambda(k) (t-s)} \hat{\bs{F}}(k,s) \diff s \diff k
  \\
  - \frac{1}{2\pi} \int_{-\infty}^\infty \e^{\im k x} S(k) \bs{a}(k,t) \diff k
  + \frac{1}{2\pi} \int_{-\infty}^\infty \e^{\im k x} \e^{-\im k l} S(k) \bs{b}(k,t) \diff k.
\end{multline}
This formula is not yet the final representation of the solution to the problem~\eqref{eq:IBVPSystem}--\eqref{eq:IBVPBoundCond}. Indeed, the formula contains the vectors $\bs{a}$ and $\bs{b}$. The vectors depend on the functions $\tilde{e}_0^{(1,2)}$ and $\tilde{q}_{0,1}^{(1,2)}$ (see Eq.~\eqref{eq:ab}) that are time transforms of the boundary values of the unknown functions $e$ and $q$ and their spatial derivatives (see Eq.~\eqref{eq:tildeeq}). Depending on the coefficients in the boundary conditions~\eqref{eq:IBVPBoundCond}, some of these boundary functions or all of them are unknown and need to be found.

\subsection{Deformation of the integration path}

In order to obtain a genuine representation of the solution, we first need to deform the integration path in the integrals involving the vectors  $\bs{a}$ and $\bs{b}$ in Eq.~\eqref{eq:SolutionRaw} from the real axis to appropriate paths in the complex $k$-plane.
We define the domains $D^\pm = D \cap \CC^\pm$, where $\CC^\pm = \{ k \in \CC : \Imag k \gtrless 0 \}$, and
\begin{equation*}
  D
  =
  \{ k \in \CC : \Real \omega_1^{}(k) < 0 \text{ or } \Real \omega_2^{}(k) < 0 \}.
\end{equation*}
Using the same reasoning as in Refs.~\cite{Fokas:2002IMAJAM, FokasPelloni:2005IMAJAM, Fokas:2008, DeconinckEtAl:2014SIAMRev, DeconinckEtAl:2018QAM}, we transform the formula~\eqref{eq:SolutionRaw} into the formula
\begin{multline}
  \label{eq:SolutionRawDeformedPath}
  \bs{u}(x,t)
  =
  \frac{1}{2\pi} \int_{-\infty}^\infty \e^{\im k x} \e^{-\varLambda(k) t} \hat{\bs{u}}_0^{}(k) \diff k
  + \frac{1}{2\pi} \int_{-\infty}^\infty \e^{\im k x} \int_0^t \e^{-\varLambda(k) (t-s)} \hat{\bs{F}}(k,s) \diff s \diff k
  \\
  - \frac{1}{2\pi} \int_{\pd D_0^{+}}^{} \e^{\im k x} S(k) \bs{a}(k,t) \diff k
  - \frac{1}{2\pi} \int_{\pd D_0^{-}}^{} \e^{-\im k (l - x)} S(k) \bs{b}(k,t) \diff k,
\end{multline}
where the boundaries $\pd D^{\pm}$ of the domains $D^{\pm}$ are given by
\begin{equation}
  \label{eq:pdD}
  \pd D^{\pm}
  =
  \left\{ k = k_{\mrm{R}}^{} + \im k_{\mrm{I}}^{} : k_{\mrm{I}}^{} = \pm \sqrt{\dfrac{\big( \mu^2 |k|^2 - 1 \big)^2 + 4 \mu^2 k_{\mrm{R}}^2}{\big( \mu^2 |k|^2 - 1 \big)^2 + 4C \mu^2 k_{\mrm{R}}^2}} |k_{\mrm{R}}^{}|, \ k_{\mrm{R}}^{} \in (-\infty, \infty) \right\},
\end{equation}
$C = \alpha \tau / \mu^2$ (see Appendix~\ref{sec:PathsOfIntegration}), and we assume that the path of integration along the boundary of a domain is directed in such a way that the domain is to the left of the boundary when it is traversed. Besides, we deform the paths of integration from $\pd D^{\pm}$ to $\pd D_0^{\pm}$ in the vicinity of the origin.
Note that if $C=1$, then the boundaries are given by
\begin{equation*}
  \pd D^{\pm}
  =
  \left\{ k = k_{\mrm{R}}^{} + \im k_{\mrm{I}}^{} : k_{\mrm{I}}^{} = \pm |k_{\mrm{R}}^{}|, \ k_{\mrm{R}}^{} \in (-\infty, \infty) \right\}.
\end{equation*}
Eq.~\eqref{eq:pdD} implies also that
\begin{equation*}
  \left| \frac{k_{\mrm{I}}^{}}{k_{\mrm{R}}^{}} \right| \to 1
  \quad\text{as}\quad
  k_{\mrm{R}}^{} \to 0
  \enskip\text{or}\enskip
  k_{\mrm{R}}^{} \to \infty,
\end{equation*}
see Appendix~\ref{sec:PathsOfIntegration}.

\subsection{Using the symmetry of the dispersion relations}

Next, we need to use the symmetry of the dispersion relations.
They are invariant under the transformation $k \to -k$: $\omega_{1,2}^{}(-k) = \omega_{1,2}^{}(k)$.
Using this transformation in the global relation~\eqref{eq:GlobRelM} we obtain another global relation (note that $\varOmega(-k) = \varOmega(k)$)
\begin{multline}
  \label{eq:GlobRelMMinusK}
  S^{-1}(-k) \hat{\bs{u}}(-k,t)
  - \e^{-\varOmega(k) t} S^{-1}(-k) \hat{\bs{u}}_0^{}(-k)
  + \bs{a}(-k,t)
  - \e^{\im k l} \bs{b}(-k,t)
  \\
  =
  \tilde{\bs{F}}(-k,t),
  \quad
  k \in \CC,
\end{multline}
To make formulae more compact, we define the vectors
\begin{equation}
  \label{eq:U}
  \bs{U}(k,t)
  \equiv
  \begin{pmatrix}
    U_1^{}(k,t) \\[2pt] U_2^{}(k,t)
  \end{pmatrix}
  \defin
  S^{-1}(k) \hat{\bs{u}}(k,t)
  =
  \begin{pmatrix}
    -\omega_2^{}
    \\[2pt]
    \omega_1^{}
  \end{pmatrix}
  \hat{e}(k,t)
  +
  \begin{pmatrix}
    \im k
    \\[2pt]
    -\im k
  \end{pmatrix}
  \hat{q}(k,t),
\end{equation}
\begin{equation}
  \label{eq:varPhi}
  \bs\varPhi(k)
  \equiv
  \begin{pmatrix}
    \varPhi_1^{}(k) \\[2pt] \varPhi_2^{}(k)
  \end{pmatrix}
  \defin
  \e^{-\varOmega(k) t} S^{-1}(k) \hat{\bs{u}}_0^{}(k)
  =
  \begin{pmatrix}
    -\omega_2^{} \e^{-\omega_1^{} t}
    \\[2pt]
    \omega_1^{} \e^{-\omega_2^{} t}
  \end{pmatrix}
  \hat{\varphi}(k)
  +
  \begin{pmatrix}
    \im k \e^{-\omega_1^{} t}
    \\[2pt]
    - \im k \e^{-\omega_2^{} t}
  \end{pmatrix}
  \hat{\psi}(k).
\end{equation}
Then the global relations~\eqref{eq:GlobRelM}, \eqref{eq:GlobRelMMinusK} can be written in the form
\begin{subequations}
  \label{eq:GlobRelUV}
  \begin{align}
    \label{eq:GlobRelMUV}
    \bs{a}(k,t)
    - \e^{-\im k l} \bs{b}(k,t)
    & =
    - \bs{U}(k,t)
    + \bs\varPhi(k,t)
    + \tilde{\bs{F}}(k,t),
    \quad
    k \in \CC,
    \\
    \label{eq:GlobRelMMinusKUV}
    \bs{a}(-k,t)
    - \e^{\im k l} \bs{b}(-k,t)
    & =
    - \bs{U}(-k,t)
    + \bs\varPhi(-k,t)
    + \tilde{\bs{F}}(-k,t),
    \quad
    k \in \CC,
  \end{align}
\end{subequations}
respectively, where $\bs{a}$ and $\bs{b}$ are given by Eq.~\eqref{eq:ab}, and $\tilde{\bs{F}}$ is given by Eq.~\eqref{eq:tildeF}.

Now we apply the time transform to the boundary conditions~\eqref{eq:IBVPBoundCond} and obtain the relations
\begin{equation}
  \label{eq:IBVPBoundCondTilde}
  \gamma_0^{} \tilde{e}_0^{(m)}(0,t)
  + \delta_0^{} \tilde{q}_0^{(m)}(0,t)
  =
  \tilde{g}^{(m)}(t),
  \quad
  \gamma_l^{} \tilde{e}_0^{(m)}(l,t)
  + \delta_l^{} \tilde{q}_0^{(m)}(l,t)
  =
  \tilde{h}^{(m)}(t),
\end{equation}
where
\begin{subequations}
  \label{eq:tildegh}
  \begin{align}
    \label{eq:tildeg}
    \tilde{g}^{(m)}(t)
    \equiv
    \tilde{g}(\omega_m^{}(k),t)
    =
    \int_0^t \e^{-\omega_m^{}(k) (t-s)} g(s) \diff s,
    \\[2pt]
    \label{eq:tildeh}
    \tilde{h}^{(m)}(t)
    \equiv
    \tilde{h}(\omega_m^{}(k),t)
    =
    \int_0^t \e^{-\omega_m^{}(k) (t-s)} h(s) \diff s.
  \end{align}
\end{subequations}
Note that $\tilde{g}^{(m)}$ and $\tilde{h}^{(m)}$ do not change under the transformation $k \to -k$.

To avoid too cumbersome general calculations, below we will limit ourselves to the most interesting boundary conditions from a physical point of view. We will limit ourselves to boundary conditions that determine the heat flux or express Newton's law. We will not consider the conditions determining the boundary temperature. Therefore, we assume that $\delta_0^{} = 1$, $\delta_l^{} = -1$ and  $\gamma_{0,l}^{} \geq 0$. In this case the boundary conditions take the form
\begin{equation}
  \label{eq:IBVPBoundCondM}
  \left. (\gamma_0^{} e + q) \right|_{x=0}^{}
  =
  g(t),
  \quad
  \left. (\gamma_l^{} e - q) \right|_{x=l}^{}
  =
  h(t),
\end{equation}
and the relations~\eqref{eq:IBVPBoundCondTilde} take the form
\begin{equation}
  \label{eq:IBVPBoundCondMTilde}
    \tilde{q}_0^{(m)}(0,t)
    =
    \tilde{g}^{(m)}(t)
    - \gamma_0^{} \tilde{e}_0^{(m)}(0,t),
    \quad
    \tilde{q}_0^{(m)}(l,t)
    =
    \gamma_l^{} \tilde{e}_0^{(m)}(l,t)
    - \tilde{h}^{(m)}(t).
\end{equation}
Besides, Eqs.~\eqref{eq:IBVPEqA} and~\eqref{eq:IBVPInitCond} imply the relations
\begin{equation}
  \label{eq:tildeqone}
  \tilde{q}_1^{(m)}(x,t)
  =
  \omega_m^{} \tilde{e}_0^{(m)}(x,t)
  - e(x,t)
  + \varphi(x) \e^{-\omega_m^{} t}
  + \tilde{f}^{(m)}(x,t),
\end{equation}
where
\begin{equation}
  \label{eq:tildef}
  \tilde{f}^{(m)}(x,t)
  \equiv
  \tilde{f}(\omega_m^{}(k),x,t)
  =
  \int_0^t \e^{-\omega_m^{}(k) (t-s)} f(x,s) \diff s
\end{equation}
is the time transform of $f$.

Using the relations~\eqref{eq:IBVPBoundCondMTilde} and~\eqref{eq:tildeqone} we can express the vectors $\bs{a}$ and $\bs{b}$ (Eq.~\eqref{eq:ab}) via $\tilde{e}_0^{(1,2)}$, $\tilde{g}^{(m)}$ and $\tilde{h}^{(m)}$:
\begin{subequations}
  \label{eq:abMatrixM}
  \begin{align}
    \label{eq:aMatrixM}
    \bs{a}(k,t)
    =&
       \begin{pmatrix}
         - \sigma_0^{(1)}(k) \tilde{e}_0^{(1)}(0,t)
         + (\omega_2^{} - \theta k^2) \tilde{g}^{(1)}(t)
         \\[3pt]
         \sigma_0^{(2)}(k) \tilde{e}_0^{(2)}(0,t)
         - (\omega_1^{} - \theta k^2) \tilde{g}^{(2)}(t)
       \end{pmatrix}
    + \im k \theta
    \begin{pmatrix}
      - e(0,t)
      + L_1^{}(k,t)
      \\[3pt]
      e(0,t)
      - L_2^{}(k,t)
    \end{pmatrix},
    \\[4pt]
    \label{eq:bMatrixM}
    \bs{b}(k,t)
    =&
       \begin{pmatrix}
         \sigma_l^{(1)}(-k) \tilde{e}_0^{(1)}(l,t)
         - (\omega_2^{} - \theta k^2) \tilde{h}^{(1)}(t)
         \\[3pt]
         - \sigma_l^{(2)}(-k) \tilde{e}_0^{(2)}(l,t)
         + (\omega_1^{} - \theta k^2) \tilde{h}^{(2)}(t)
       \end{pmatrix}
    + \im k \theta
    \begin{pmatrix}
      - e(l,t)
      + R_1^{}(k,t)
      \\[3pt]
      e(l,t)
      - R_2^{}(k,t)
    \end{pmatrix},
  \end{align}
\end{subequations}
where
\begin{equation}
  \label{eq:sigma}
  \begin{split}
  \sigma_{0,l}^{(1)}(k)
  =
  \im k (\beta - \theta \omega_1^{})
  + \gamma_{0,l}^{} (\omega_2^{} - \theta k^2),
  \\
  \sigma_{0,l}^{(2)}(k)
  =
  \im k (\beta - \theta \omega_2^{})
  + \gamma_{0,l}^{} (\omega_1^{} - \theta k^2),
  \end{split}
\end{equation}
\begin{equation}
  \label{eq:LR}
  \begin{split}
    L_1^{}(k,t)
    =
    \varphi(0) \e^{-\omega_1^{} t}
    + \tilde{f}^{(1)}(0,t),
    \quad
    R_1^{}(k,t)
    =
    \varphi(l) \e^{-\omega_1^{} t}
    + \tilde{f}^{(1)}(l,t),
    \\
    L_2^{}(k,t)
    =
    \varphi(0) \e^{-\omega_2^{} t}
    + \tilde{f}^{(2)}(0,t),
    \quad
    R_2^{}(k,t)
    =
    \varphi(l) \e^{-\omega_2^{} t}
    + \tilde{f}^{(2)}(l,t).
  \end{split}
\end{equation}
The vectors $\bs{a}$ and $\bs{b}$ depend on the functions $\tilde{e}_0^{(1,2)}$, which are still unknown. To find them, we use the system~\eqref{eq:GlobRelUV} of the global relations and the representations~\eqref{eq:abMatrixM}.
The system~\eqref{eq:GlobRelUV}, consisting of four equations, can be rewritten as the system of two algebraic systems
\begin{subequations}
  \label{eq:SystemsForE}
  \begin{align}
    \label{eq:SystemsForEOne}
    &
    \begin{pmatrix}
      \sigma_0^{(1)}(k)
      & \e^{-\im k l} \sigma_l^{(1)}(-k)
      \\[4pt]
      \sigma_0^{(1)}(-k)
      & \e^{\im k l} \sigma_l^{(1)}(k)
    \end{pmatrix}
    \begin{pmatrix}
      \tilde{e}_0^{(1)}(0,t)
      \\[4pt]
      \tilde{e}_0^{(1)}(l,t)
    \end{pmatrix}
    \\
    \notag
    &\hspace{90pt}=
    \begin{pmatrix}
      U_1^{}(k,t) - V_1^{}(k,t) - W_1^{}(k,t) - E(k,t)
      \\[4pt]
      U_1^{}(-k,t) - V_1^{}(-k,t) - W_1^{}(-k,t) - E(-k,t)
    \end{pmatrix},
    \\[2pt]
    \label{eq:SystemsForETwo}
    &
    \begin{pmatrix}
      \sigma_0^{(2)}(k)
      & \e^{-\im k l} \sigma_l^{(2)}(-k)
      \\[2pt]
      \sigma_0^{(2)}(-k)
      & \e^{\im k l} \sigma_l^{(2)}(k)
    \end{pmatrix}
    \begin{pmatrix}
      \tilde{e}_0^{(2)}(0,t)
      \\[2pt]
      \tilde{e}_0^{(2)}(l,t)
    \end{pmatrix}
    \\
    \notag
    &\hspace{90pt}=
    -
    \begin{pmatrix}
      U_2^{}(k,t) - V_2^{}(k,t) - W_2^{}(k,t) - E(k,t)
      \\[2pt]
      U_2^{}(-k,t) - V_2^{}(-k,t) - W_2^{}(-k,t) - E(-k,t)
    \end{pmatrix}
  \end{align}
\end{subequations}
(the equations of the system~\eqref{eq:SystemsForEOne} are the first equations (raws) of the systems~\eqref{eq:GlobRelMUV} and~\eqref{eq:GlobRelMMinusKUV}, and the equations of the system~\eqref{eq:SystemsForETwo} are the second equations (raws) of these systems), where
\begin{equation}
  \label{eq:V}
  \begin{split}
    V_1^{}(k,t)
    &=
    \varPhi_1^{}(k) + \tilde{F}_1^{}(k,t)
    - (\omega_2^{} - \theta k^2) \big[ \tilde{g}^{(1)}(t) + \e^{-\im k l} \tilde{h}^{(1)}(t) \big],
    \\
    V_2^{}(k,t)
    &=
    \varPhi_2^{}(k) + \tilde{F}_2^{}(k,t)
    + (\omega_1^{} - \theta k^2) \big[ \tilde{g}^{(2)}(t) + \e^{-\im k l} \tilde{h}^{(2)}(t) \big],
  \end{split}
\end{equation}
$\varPhi_{1,2}^{}$ are given by Eq.~\eqref{eq:varPhi}, $\tilde{F}_{1,2}^{}$ are given by Eq.~\eqref{eq:tildeF}, $\tilde{g}^{(1,2)}$ and $\tilde{h}^{(1,2)}$ are given by Eq.~\eqref{eq:tildegh},
\begin{equation}
  \label{eq:W}
  \begin{split}
    W_1^{}(k,t)
    &=
    - \im k \theta \left[
      L_1^{}(k,t)
      - \e^{-\im k l} R_1^{}(k,t)
    \right]
    \\
    \notag
    &\equiv
    - \im k \theta \left[
      \varphi(0) \e^{-\omega_1^{} t}
      + \tilde{f}^{(1)}(0,t)
      - \e^{-\im k l} \left(
        \varphi(l) \e^{-\omega_1^{} t}
        + \tilde{f}^{(1)}(l,t)
      \right)
    \right],
    \\
    W_2^{}(k,t)
    &=
    \im k \theta \left[
      L_2^{}(k,t)
      - \e^{-\im k l} R_2^{}(k,t)
    \right]
    \\
    \notag
    &\equiv
    \im k \theta \left[
      \varphi(0) \e^{-\omega_2^{} t}
      + \tilde{f}^{(2)}(0,t)
      - \e^{-\im k l} \left(
        \varphi(l) \e^{-\omega_2^{} t}
        + \tilde{f}^{(2)}(l,t)
      \right)
    \right],
  \end{split}
\end{equation}
and
\begin{equation}
  \label{eq:E}
  E(k,t)
  =
  \im k \theta \big[ e(0,t) - \e^{-\im k l} e(l,t) \big].
\end{equation}
The inverse of the matrices in the systems~\eqref{eq:SystemsForEOne} and~\eqref{eq:SystemsForETwo} are given by
\begin{equation*}
  \frac{1}{\varDelta_1^{}(k)}
  \begin{pmatrix}
    \e^{\im k l} \sigma_l^{(1)}(k)
    & - \e^{-\im k l} \sigma_l^{(1)}(-k)
    \\[4pt]
    - \sigma_0^{(1)}(-k)
    & \sigma_0^{(1)}(k)
  \end{pmatrix}
  \quad\text{and}\quad
  \frac{1}{\varDelta_2^{}(k)}
  \begin{pmatrix}
    \e^{\im k l} \sigma_l^{(2)}(k)
    & - \e^{-\im k l} \sigma_l^{(2)}(-k)
    \\[4pt]
    - \sigma_0^{(2)}(-k)
    & \sigma_0^{(2)}(k)
  \end{pmatrix},
\end{equation*}
respectively, where
\begin{equation}
  \label{eq:DeltaOneTwo}
  \begin{split}
    \varDelta_1^{}(k)
    &=
    \e^{\im k l} \sigma_0^{(1)}(k) \sigma_l^{(1)}(k)
    - \e^{-\im k l} \sigma_0^{(1)}(-k) \sigma_l^{(1)}(-k),
    \\[2pt]
    \varDelta_2^{}(k)
    &=
    \e^{\im k l} \sigma_0^{(2)}(k) \sigma_l^{(2)}(k)
    - \e^{-\im k l} \sigma_0^{(2)}(-k) \sigma_l^{(2)}(-k).
  \end{split}
\end{equation}
Solving the systems~\eqref{eq:SystemsForE}, we find the functions $\tilde{e}_0^{(1,2)}$:
\begin{subequations}
  \label{eq:ee}
  \begin{align}
    \label{eq:eOne}
    \begin{pmatrix}
      \tilde{e}_0^{(1)}(0,t)
      \\[2pt]
      \tilde{e}_0^{(1)}(l,t)
    \end{pmatrix}
    =
    &
    \frac{1}{\varDelta_1^{}(k)} \Bigg[
    \begin{pmatrix}
      \e^{\im k l} \sigma_l^{(1)}(k)
      \\[2pt]
      - \sigma_0^{(1)}(-k)
    \end{pmatrix}
    \left[ U_1^{}(k,t) - V_1^{}(k,t) - W_1^{}(k,t) - E(k,t) \right]
    \\[2pt]
    \notag
    &+
    \begin{pmatrix}
      - \e^{-\im k l} \sigma_l^{(1)}(-k)
      \\[2pt]
      \sigma_0^{(1)}(k)
    \end{pmatrix}
    \left[ U_1^{}(-k,t) - V_1^{}(-k,t) - W_1^{}(-k,t) - E(-k,t) \right]
    \Bigg],
    \\[4pt]
    \label{eq:eTwo}
    \begin{pmatrix}
      \tilde{e}_0^{(2)}(0,t)
      \\[2pt]
      \tilde{e}_0^{(2)}(l,t)
    \end{pmatrix}
    =
    &
    \frac{1}{\varDelta_2^{}(k)} \Bigg[
    \begin{pmatrix}
      - \e^{\im k l} \sigma_l^{(2)}(k)
      \\[2pt]
      \sigma_0^{(2)}(-k)
    \end{pmatrix}
    \left[ U_2^{}(k,t) - V_2^{}(k,t) - W_2^{}(k,t) - E(k,t) \right]
    \\[2pt]
    \notag
    &+
    \begin{pmatrix}
      \e^{-\im k l} \sigma_l^{(2)}(-k)
      \\[2pt]
      - \sigma_0^{(2)}(k)
    \end{pmatrix}
    \left[ U_2^{}(-k,t) - V_2^{}(-k,t) - W_2^{}(-k,t) - E(-k,t) \right]
    \Bigg],
  \end{align}
\end{subequations}
where $U_{1,2}^{}$ are given by Eq.~\eqref{eq:U}, $V_{1,2}^{}$ are given by Eq.~\eqref{eq:V}, $W_{1,2}^{}$ are given by Eq.~\eqref{eq:W} and $E$ is given by Eq.~\eqref{eq:E}.
Note that the functions $U_{1,2}^{}$ and $E$ are still unknown, while the functions $V_{1,2}$ and $W_{1,2}$ are known, they contain only known functions.

\subsection{The solution}
\label{sec:Solution}

The formula~\eqref{eq:SolutionRawDeformedPath} is also not the final representation of the solution, because the vectors $\bs{a}$ and $\bs{b}$ depend on the unknown functions $e$ and $q$ through the functions $U_{1,2}$ and $E$ (see Eqs.~\eqref{eq:abMatrixM}, \eqref{eq:ee}, \eqref{eq:U} and \eqref{eq:E}). However, it can be shown that integrals containing the functions $\hat{e}(k,t)$, $\hat{q}(k,t)$, $e(0,t)$ and $e(l,t)$ are equal to zero, see, \eg, Refs.~\cite{FokasPelloni:2005IMAJAM, Fokas:2008, DeconinckEtAl:2014SIAMRev, DeconinckEtAl:2018QAM}, and, therefore, the final formula does not contain these integrals.

As a result, the solution is given by
\begin{multline}
  \label{eq:Solution}
  \bs{u}(x,t)
  =
  \frac{1}{2\pi} \int_{-\infty}^\infty \e^{\im k x} \e^{-\varLambda(k) t} \hat{\bs{u}}_0^{}(k) \diff k
  + \frac{1}{2\pi} \int_{-\infty}^\infty \e^{\im k x} \int_0^t \e^{-\varLambda(k) (t-s)} \hat{\bs{F}}(k,s) \diff s \diff k
  \\
  - \frac{1}{2\pi} \int_{\pd D_0^{+}}^{} \e^{\im k x} S(k) \bs{a}'(k,t) \diff k
  - \frac{1}{2\pi} \int_{\pd D_0^{-}}^{} \e^{\im k x} \e^{-\im k l} S(k) \bs{b}'(k,t) \diff k,
\end{multline}
where $\hat{\bs{u}}_0^{}(k)$ is given by Eq.~\eqref{eq:HatUZero}, $\hat{\bs{F}}(k,t)$ is given by Eq.~\eqref{eq:HatF}, the vectors $\bs{a}'$ and $\bs{b}'$ are given by
\begin{subequations}
  \label{eq:abPrimeMatrixM}
  \begin{align}
    \label{eq:aPrimeMatrixM}
    \bs{a}'(k,t)
    =
    \begin{pmatrix}
      - \sigma_0^{(1)}(k) \tilde{e}_0^{(1)}(0,t)
      + (\omega_2^{} - \theta k^2) \tilde{g}^{(1)}(t)
      \\[3pt]
      \sigma_0^{(2)}(k) \tilde{e}_0^{(2)}(0,t)
      - (\omega_1^{} - \theta k^2) \tilde{g}^{(2)}(t)
    \end{pmatrix}
    + \im k \theta
    \begin{pmatrix}
      L_1^{}(k,t)
      \\[3pt]
      - L_2^{}(k,t)
    \end{pmatrix},
    \\[4pt]
    \label{eq:bPrimeMatrixM}
    \bs{b}'(k,t)
    =
    \begin{pmatrix}
      \sigma_l^{(1)}(-k) \tilde{e}_0^{(1)}(l,t)
      - (\omega_2^{} - \theta k^2) \tilde{h}^{(1)}(t)
      \\[3pt]
      - \sigma_l^{(2)}(-k) \tilde{e}_0^{(2)}(l,t)
      + (\omega_1^{} - \theta k^2) \tilde{h}^{(2)}(t)
    \end{pmatrix}
    + \im k \theta
    \begin{pmatrix}
      R_1^{}(k,t)
      \\[3pt]
      - R_2^{}(k,t)
    \end{pmatrix}
  \end{align}
\end{subequations}
(these are Eqs.~\eqref{eq:abMatrixM} with $e(0,t) = 0$ and $e(l,t) = 0$ ),
the functions $\tilde{e}_0^{(1,2)}$ are given by
\begin{subequations}
  \label{eq:eePrime}
  \begin{align}
    \label{eq:eOnePrime}
    \begin{pmatrix}
      \tilde{e}_0^{(1)}(0,t)
      \\[2pt]
      \tilde{e}_0^{(1)}(l,t)
    \end{pmatrix}
    =
    \frac{1}{\varDelta_1^{}(k)} \Bigg[
    &
    \begin{pmatrix}
      -\e^{\im k l} \sigma_l^{(1)}(k)
      \\[2pt]
      \sigma_0^{(1)}(-k)
    \end{pmatrix}
    \left[ V_1^{}(k,t) + W_1^{}(k,t) \right]
    \\[2pt]
    \notag
    &+
    \begin{pmatrix}
      \e^{-\im k l} \sigma_l^{(1)}(-k)
      \\[2pt]
      -\sigma_0^{(1)}(k)
    \end{pmatrix}
    \left[ V_1^{}(-k,t) + W_1^{}(-k,t) \right]
    \Bigg],
    \\[4pt]
    \label{eq:eTwoPrime}
    \begin{pmatrix}
      \tilde{e}_0^{(2)}(0,t)
      \\[2pt]
      \tilde{e}_0^{(2)}(l,t)
    \end{pmatrix}
    =
    \frac{1}{\varDelta_2^{}(k)} \Bigg[
    &
    \begin{pmatrix}
      \e^{\im k l} \sigma_l^{(2)}(k)
      \\[2pt]
      -\sigma_0^{(2)}(-k)
    \end{pmatrix}
    \left[ V_2^{}(k,t) + W_2^{}(k,t) \right]
    \\[2pt]
    \notag
    &+
    \begin{pmatrix}
      -\e^{-\im k l} \sigma_l^{(2)}(-k)
      \\[2pt]
      \sigma_0^{(2)}(k)
    \end{pmatrix}
    \left[ V_2^{}(-k,t) + W_2^{}(-k,t) \right]
    \Bigg]
  \end{align}
\end{subequations}
(these are Eqs.~\eqref{eq:ee} with $U_{1,2}^{} = 0$ and $E = 0$).

\section{A numerical example}
\label{sec:Example}

To illustrate the representation~\eqref{eq:Solution} of the solution, we consider here a particular case, in which heat sources are absent ($f = 0$), the initial state is zero ($\varphi = 0$, $\psi = 0$), the heat flux is prescribed at the left boundary ($\gamma_0^{} = 0$), and heat exchange at the right boundary obeys Newton's law, while the ambient temperature is zero ($h = 0$),
In this case both equations are homogeneous:
\begin{subequations}
  \label{eq:IBVPEqsCaseOne}
  \begin{align}
    &\pd_t^{} e + \pd_x^{} q
    = 0,
    \quad
    x \in (0,l),
    \quad
    t > 0,
    \\
    &\tau \pd_t^{} q
    + q
      - \mu^2 \pd_x^2 q
    + \alpha \pd_x^{} e
    = 0,
    \quad
    x \in (0,l),
    \quad
    t > 0.
  \end{align}
\end{subequations}
the initial conditions are also homogeneous:
\begin{equation}
  \label{eq:IBVPInitCondCaseOne}
  \left. e \right|_{t=0}^{}
  = 0,
  \quad
  \left. q \right|_{t=0}^{}
  = 0,
\end{equation}
and the boundary conditions take the form
\begin{equation}
  \label{eq:IBVPBoundCondCaseOne}
  \left. q \right|_{x=0}^{}
  =
  g(t),
  \quad
  \left. q \right|_{x=l}^{}
  =
  \left. \gamma_l^{} e \right|_{x=l}^{}
\end{equation}
(the second boundary condition expresses Newton's law, when the ambient temperature is zero).

The solution to the problem~\eqref{eq:IBVPEqsCaseOne}--\eqref{eq:IBVPBoundCondCaseOne} is given by
\begin{equation}
  \label{eq:SolutionCaseOne}
  \bs{u}(x,t)
  =
  - \frac{1}{2\pi} \int_{\pd D_0^{+}}^{} \e^{\im k x} S(k) \bs{a}'(k,t) \diff k
  - \frac{1}{2\pi} \int_{\pd D_0^{-}}^{} \e^{\im k x} \e^{-\im k l} S(k) \bs{b}'(k,t) \diff k,
\end{equation}
where the vectors $\bs{a}'$ and $\bs{b}'$ are given by
\begin{equation*}
  \begin{aligned}
    \bs{a}'(k,t)
    &=
    \begin{pmatrix}
      - \im k (\beta - \theta \omega_1^{}) \tilde{e}_0^{(1)}(0,t)
      + (\omega_2^{} - \theta k^2) \tilde{g}^{(1)}(t)
      \\[3pt]
      \im k (\beta - \theta \omega_2^{}) \tilde{e}_0^{(2)}(0,t)
      - (\omega_1^{} - \theta k^2) \tilde{g}^{(2)}(t)
    \end{pmatrix},
    \\[4pt]
    \bs{b}'(k,t)
    &=
    \begin{pmatrix}
      \sigma_l^{(1)}(-k) \tilde{e}_0^{(1)}(l,t)
      \\[3pt]
      - \sigma_l^{(2)}(-k) \tilde{e}_0^{(2)}(l,t)
    \end{pmatrix},
  \end{aligned}
\end{equation*}
(see Eq.~\eqref{eq:abPrimeMatrixM}),
the coefficients $\sigma_l^{(1,2)}(k)$ are given by Eq.~\eqref{eq:sigma},
the functions $\tilde{e}_0^{(1,2)}$ are given by
\begin{subequations}
  \label{eq:eePrimeCaseOne}
  \begin{align}
    \label{eq:eOnePrimeCaseOne}
    \begin{pmatrix}
      \im k (\beta - \theta \omega_1^{}) \tilde{e}_0^{(1)}(0,t)
      \\[3pt]
      \tilde{e}_0^{(1)}(l,t)
    \end{pmatrix}
    &=
    \frac{1}{\varDelta_1'(k)}
    \begin{pmatrix}
      \e^{\im k l} \sigma_l^{(1)}(k)
      - \e^{-\im k l} \sigma_l^{(1)}(-k)
      \\[3pt]
      2
    \end{pmatrix}
    (\omega_2^{} - \theta k^2) \tilde{g}^{(1)}(t),
    \\[4pt]
    \label{eq:eTwoPrimeCaseOne}
    \begin{pmatrix}
      \im k (\beta - \theta \omega_2^{}) \tilde{e}_0^{(2)}(0,t)
      \\[3pt]
      \tilde{e}_0^{(2)}(l,t)
    \end{pmatrix}
    &=
    \frac{1}{\varDelta_2'(k)}
    \begin{pmatrix}
      \e^{\im k l} \sigma_l^{(2)}(k)
      - \e^{-\im k l} \sigma_l^{(2)}(-k)
      \\[3pt]
      2
    \end{pmatrix}
    (\omega_1^{} - \theta k^2) \tilde{g}^{(2)}(t),
  \end{align}
\end{subequations}
the coefficients $\varDelta_{1,2}'(k)$ are given by
\begin{equation*}
  \begin{split}
    \varDelta_1'(k)
    =
    \e^{\im k l} \sigma_l^{(1)}(k)
    + \e^{-\im k l} \sigma_l^{(1)}(-k),
    \\
    \varDelta_2'(k)
    =
    \e^{\im k l} \sigma_l^{(2)}(k)
    + \e^{-\im k l} \sigma_l^{(2)}(-k).
  \end{split}
\end{equation*}
Using the relations~\eqref{eq:eePrimeCaseOne}, we see that
\begin{equation*}
  \begin{split}
    - \im k (\beta - \theta \omega_1^{}) \tilde{e}_0^{(1)}(0,t)
    + (\omega_2^{} - \theta k^2) \tilde{g}^{(1)}(t)
    =
    2 \e^{-\im k l} \frac{\sigma_l^{(1)}(-k)}{\varDelta_1'(k)} (\omega_2^{} - \theta k^2) \tilde{g}^{(1)}(t),
    \\
    - \im k (\beta - \theta \omega_2^{}) \tilde{e}_0^{(2)}(0,t)
    + (\omega_1^{} - \theta k^2) \tilde{g}^{(2)}(t)
    =
    2 \e^{-\im k l} \frac{\sigma_l^{(2)}(-k)}{\varDelta_2'(k)} (\omega_1^{} - \theta k^2) \tilde{g}^{(2)}(t).
  \end{split}
\end{equation*}
Hence
\begin{multline}
  \label{eq:SaSbPrimeCaseOneM}
  S(k) \bs{a}'(k,t)
  =
  \e^{-\im k l} S(k) \bs{b}'(k,t)
  \\[2pt]
  =
  \frac{2 \e^{-\im k l}}{\omega_1^{} - \omega_2^{}} \bigg[
  \begin{pmatrix}
    1
    \\
    \omega_1^{} / \im k
  \end{pmatrix}
  \frac{\sigma_l^{(1)}(-k)}{\varDelta_1'(k)} (\omega_2^{} - \theta k^2) \tilde{g}^{(1)}(t)
  -
  \begin{pmatrix}
    1
    \\
    \omega_2^{} / \im k
  \end{pmatrix}
  \frac{\sigma_l^{(2)}(-k)}{\varDelta_2'(k)} (\omega_1^{} - \theta k^2) \tilde{g}^{(2)}(t)
  \bigg].
\end{multline}
As a result, we conclude that the solution to the problem~\eqref{eq:IBVPEqsCaseOne}--\eqref{eq:IBVPBoundCondCaseOne} is given by the formula~\eqref{eq:SolutionCaseOne}, where $S \bs{a}'$ and $S \bs{b}'$ are given by Eq.~\eqref{eq:SaSbPrimeCaseOneM}.

As a function, determining the heat flux on the left boundary, we choose the function
\begin{equation}
  \label{eq:g}
  g(t)
  =
  \frac{2}{\tau_\Delta^{}} \sin^2 \!\left( \frac{\pi t}{\tau_\Delta^{}} \right) \mathbbm{1}_{(0,\tau_\Delta^{})}^{}(t),
\end{equation}
where
\begin{equation*}
  \label{eq:Indicator}
  \mathbbm{1}_A^{}(t)
  =
  \begin{cases}
    1, & t \in A,
    \\
    0, & t \notin A,
  \end{cases}
\end{equation*}
is the indicator function of a set $A$. Note that $\int_0^\infty g(t) \diff t \equiv \int_0^{\tau_\Delta^{}} g(t) \diff t = 1$. The problem~\eqref{eq:IBVPEqsCaseOne}--\eqref{eq:IBVPBoundCondCaseOne} with the heat flux of the form~\eqref{eq:g} simulates laser flash experiments~\cite{Kovacs:2018IJHMT}.

\begin{figure}[!htb]
  \centering
  \includegraphics*[scale=1]{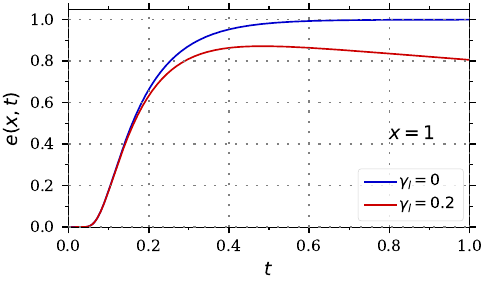}
  \caption{The solutions to the problem~\eqref{eq:IBVPEqsCaseOne}--\eqref{eq:IBVPBoundCondCaseOne} on the right boundary, \ie, for $x = 1$ with the heat flux $g(t)$ given by Eq.~\eqref{eq:g}. The values of the parameters are $\alpha = 1$, $\tau = 0.02$, $\mu^2 = 0.02$, $\gamma_l^{} = 0$ and $\gamma_l^{} = 0.2$, $l = 1$, $\tau_\Delta^{} = 0.04$.}
  \label{fig:Laser_flash_mu2=002}
\end{figure}

\begin{figure}[!htb]
  \centering
  \includegraphics*[scale=1]{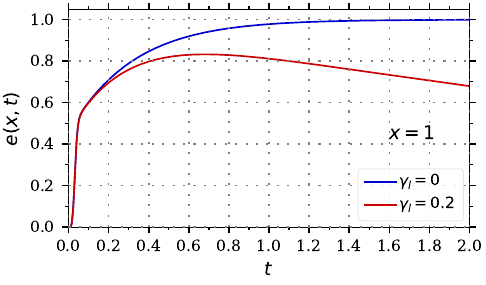}
  \caption{The solutions to the problem~\eqref{eq:IBVPEqsCaseOne}--\eqref{eq:IBVPBoundCondCaseOne} on the right boundary, \ie, for $x = 1$ with the heat flux $g(t)$ given by Eq.~\eqref{eq:g}. The values of the parameters are $\alpha = 1$, $\tau = 0.02$, $\mu^2 = 0.2$, $\gamma_l^{} = 0$ and $\gamma_l^{} = 0.2$, $l = 1$, $\tau_\Delta^{} = 0.04$.}
  \label{fig:Laser_flash_mu2=02}
\end{figure}

Figs.~\ref{fig:Laser_flash_mu2=002}, \ref{fig:Laser_flash_mu2=02} show solutions to the problem~\eqref{eq:IBVPEqsCaseOne}--\eqref{eq:IBVPBoundCondCaseOne} obtained by using the formula~\eqref{eq:SolutionCaseOne}. The values of the parameters are $\alpha = 1$, $\tau = 0.02$, $\mu^2 = 0.02$ and $\mu^2 = 0.2$, $\gamma_l^{} = 0$ and $\gamma_l^{} = 0.2$, $l = 1$, $\tau_\Delta^{} = 0.04$. Here we use dimensionless values. The solution is calculated on the right boundary, \ie, for $x = 1$. The solution to the problem with $\gamma_l^{} = 0$ (the thermally insulating boundary) and the general heat flux $g(t)$ can easily be obtained by the Fourier method with preliminary elimination of inhomogeneities in boundary conditions. This solution is given by Eq.~\eqref{eq:SolutionCaseTwoSeries} (we have checked that in this case the formulae~\eqref{eq:SolutionCaseOne} and~\eqref{eq:SolutionCaseTwoSeries} give the same numerical results). In Ref.~\cite{Kovacs:2018IJHMT} the series solution to the problem~\eqref{eq:IBVPEqsCaseOne}--\eqref{eq:IBVPBoundCondCaseOne} with  $\gamma_l^{} = 0$ was obtained for the heat flux of the form~\eqref{eq:g}.

\appendix

\section{The paths of integration in the complex plane}
\label{sec:PathsOfIntegration}

To find the paths of integration we must solve the equation
\begin{equation*}
  \Real F(k)
  = 0,
  \quad\text{where}\quad
  F(k)
  =
  1 + \mu^2 k^2 \pm \sqrt{(1 + \mu^2 k^2)^2 - 4 \alpha \tau k^2},
\end{equation*}
see Eq.~\eqref{eq:omegaOneTwo}.
Making the change of variable $\mu k = z \equiv x + \im y$ ($x = \Real z$, it is not the spatial variable, as in the main part of the paper), and introducing the parameter $C = \alpha \tau / \mu^2$, we obtain
\begin{multline*}
  F(k)
  \equiv
  G(z)
  =
  1 + z^2 \pm \sqrt{(1 + z^2)^2 - 4C z^2}
  \\
  =
  x^2 - y^2 + 1 + \im 2 xy \pm \sqrt{(x^2 - y^2 + 1 + \im 2 xy)^2 - 4C (x^2 - y^2 + \im 2 xy)}
  \\
  =
  x^2 - y^2 + 1 + \im 2 xy \pm \sqrt{\xi + \im \eta}
  \\
  =
  x^2 - y^2 + 1 + \im 2 xy \pm \sqrt{r^2 \e^{\im 2 \varphi}}
  =
  x^2 - y^2 + 1 + \im 2 xy \pm r \e^{\im \varphi},
\end{multline*}
where
\begin{equation*}
  r^2
  =
  \sqrt{\xi^2 + \eta^2},
  \quad
  \cot 2\varphi
  =
  \frac{\xi}{\eta},
  \quad
  \varphi \in \left( -\frac{\pi}{2}, \frac{\pi}{2} \right),
\end{equation*}
\begin{equation*}
  \xi
  =
  (x^2 - y^2 + 1)^2 - 4 x^2 y^2 - 4C (x^2  - y^2),
  \quad
  \eta
  =
  4 xy (x^2 - y^2 + 1 - 2C).
\end{equation*}
So the equation becomes
\begin{equation*}
  \Real G(z) = 0
  \quad\Leftrightarrow\quad
  r \cos\varphi = |x^2 - y^2 + 1|
  \quad\Leftrightarrow\quad
  r^2 \cos^2\varphi = (x^2 - y^2 + 1)^2.
\end{equation*}
Using the equality
\begin{equation*}
  \cos^2 \varphi
  =
  \frac{1}{2} \left( 1 + \sign \varphi \frac{\cot 2\varphi}{\sqrt{\cot^2 2\varphi + 1}} \right),
  \quad
  \varphi \in \left( -\frac{\pi}{2}, \frac{\pi}{2} \right),
\end{equation*}
we obtain the equation
\begin{multline*}
  \frac{r^2}{2} \left( 1 + \sign \varphi \frac{\xi}{\eta \sqrt{(\xi/\eta)^2 + 1}} \right)
  \equiv
  \frac{r^2}{2} \left( 1 + \frac{\xi}{\sqrt{\xi^2 + \eta^2}} \right)
  \equiv
  \frac{\sqrt{\xi^2 + \eta^2} + \xi}{2}
  \\
  =
  (x^2 - y^2 + 1)^2,
\end{multline*}
which is equavalent to the equation
\begin{equation*}
  \xi^2 + \eta^2
  =
  \left[ 2 (x^2 - y^2 + 1)^2 - \xi \right]^2.
\end{equation*}
As a result, we get the equation
\begin{equation}
  \label{eq:CubicEqOriginal}
  y^6
  + (x^2 - 2) y^4
  - (x^4 - 4C x^2 - 1) y^2
  - (x^6 +2 x^4 + x^2)
  = 0
\end{equation}
that is actually a cubic equation for $y^2$. This equation can be written in the form
\begin{equation*}
  (y^2 - x^2) [(x^2 + y^2 - 1)^2 + 4 x^2]
  + 4 (C - 1) x^2 y^2
  = 0,
\end{equation*}
or, equivalently
\begin{equation}
  \label{eq:yEqual}
  y
  =
  \pm \sqrt{\dfrac{(x^2 + y^2 - 1)^2 + 4 x^2}{(x^2 + y^2 - 1)^2 + 4C x^2}} |x|,
  \quad
  C = \frac{\alpha \tau}{\mu^2}.
\end{equation}
This relation expresses the dependence of $y$ on $x$ implicitly, however, it allows us to obtain qualitative information about the solutions of the equation~\eqref{eq:CubicEqOriginal}. The relation~\eqref{eq:yEqual} shows that the equation~\eqref{eq:CubicEqOriginal} has exactly two real solutions $y(x)$, and, therefore, the paths $\pd D^{\pm}$, Eq.~\eqref{eq:pdD}, are uniquely determined. If $C=1$, then the paths $\pd D^{\pm}$ are determined by the equations $k_{\mrm{I}}^{} = \pm |k_{\mrm{R}}^{}|$, respectively. Note that these paths coincide with the paths of integration for the heat equation on a finite interval, see Ref.~\cite{DeconinckEtAl:2014SIAMRev}. If $C<1$, then the path $\pd D^{+}$ passes above the path $k_{\mrm{I}}^{} = |k_{\mrm{R}}^{}|$, and the path $\pd D^{-}$ passes below the path $k_{\mrm{I}}^{} = -|k_{\mrm{R}}^{}|$. If $C>1$, then the path $\pd D^{+}$ passes below the path $k_{\mrm{I}}^{} = |k_{\mrm{R}}^{}|$ and above the real axis, and the path $\pd D^{-}$ passes above the path $k_{\mrm{I}}^{} = -|k_{\mrm{R}}^{}|$ and below the real axis.
The relation~\eqref{eq:yEqual} shows also the asymptotic behavior
\begin{equation*}
  \left| \frac{y}{x} \right| \to 1
  \quad\text{as}\quad
  x \to 0
  \enskip\text{or}\enskip
  x \to \infty.
\end{equation*}

To determine the paths $\pd D^{\pm}$ explicitly, we write the equation~\eqref{eq:CubicEqOriginal} in the form
\begin{equation}
  \label{eq:CubicEqMod}
  s^3 + a s^2 + b s + c
  =
  u^3 + P u + Q
  =
  (u - 2 A) [(u + A)^2 + B^2]
  = 0,
\end{equation}
where $y^2 \equiv s = u - a/3$,
\begin{equation*}
  a
  =
  x^2 - 2,
  \quad
  b
  =
  -(x^4 - 4C x^2 - 1),
  \quad
  c
  =
  -(x^6 +2 x^4 + x^2),
\end{equation*}
\begin{equation*}
  P
  =
  - \frac{a^2}{3}
  + b,
  \quad
  Q
  =
  2 \left( \frac{a}{3} \right)^3
  - \frac{a}{3} b
  + c,
  \quad
  A
  =
  \frac{\alpha + \beta}{2},
  \quad
  B
  =
  \sqrt{3} \,\frac{\alpha - \beta}{2},
\end{equation*}
\begin{equation*}
  \alpha
  =
  \sqrt[3]{-\frac{Q}{2} + \sqrt{D}},
  \quad
  \beta
  =
  \sqrt[3]{-\frac{Q}{2} - \sqrt{D}},
  \quad
  D
  =
  \left( \frac{P}{3} \right)^3
  + \left( \frac{Q}{2} \right)^2,
\end{equation*}
the roots $\alpha$ and $\beta$ are chosen so that the equality $\alpha \beta = - P/3$ is valid and the value $A$ is real.
The real nonnegative root of the equation~\eqref{eq:CubicEqMod} is given by
\begin{equation*}
  y^2 \equiv s
  =
  2A - \frac{a}{3}.
\end{equation*}
This relation determines the paths $\pd D^{\pm}$.

\section{Particular cases that allow series representation of the solutions}
\label{sec:ParticularCases}

\subsection{The heat sources are absent, the initial state is zero, the heat flux is prescribed on the left boundary, and the right boundary is thermally insulating}

Here we assume that heat sources are absent ($f = 0$), the initial state is zero ($\varphi = 0$, $\psi = 0$), the heat flux is prescribed on the left boundary ($\gamma_0^{} = 0$), and the right boundary is thermally insulating ($\gamma_l^{} = 0$ and $h = 0$).
In this case both equations are homogeneous:
\begin{subequations}
  \label{eq:IBVPEqsCaseTwo}
  \begin{align}
    &\pd_t^{} e + \pd_x^{} q
    = 0,
    \quad
    x \in (0,l),
    \quad
    t > 0,
    \\
    &\tau \pd_t^{} q
    + q
      - \mu^2 \pd_x^2 q
    + \alpha \pd_x^{} e
    = 0,
    \quad
    x \in (0,l),
    \quad
    t > 0.
  \end{align}
\end{subequations}
the initial conditions are also homogeneous:
\begin{equation}
  \label{eq:IBVPInitCondCaseTwo}
  \left. e \right|_{t=0}^{}
  = 0,
  \quad
  \left. q \right|_{t=0}^{}
  = 0,
\end{equation}
and the boundary conditions take the form
\begin{equation}
  \label{eq:IBVPBoundCondCaseTwo}
  \left. q \right|_{x=0}^{}
  =
  g(t),
  \quad
  \left. q \right|_{x=l}^{}
  = 0.
\end{equation}
This is a particular case of the problem~\eqref{eq:IBVPEqsCaseOne}--\eqref{eq:IBVPBoundCondCaseOne}, where $\gamma_l^{} = 0$.

In this case the solution to the problem~\eqref{eq:IBVPEqsCaseTwo}--\eqref{eq:IBVPBoundCondCaseTwo} is given by
\begin{equation}
  \label{eq:SolutionCaseTwoIntegral}
  \bs{u}(x,t)
  =
  - \frac{1}{2\pi} \int_{\pd D_0^{+}}^{} \e^{\im k x} S(k) \bs{a}'(k,t) \diff k
  - \frac{1}{2\pi} \int_{\pd D_0^{-}}^{} \e^{\im k x} \e^{-\im k l} S(k) \bs{b}'(k,t) \diff k,
\end{equation}
where
\begin{multline*}
  S(k) \bs{a}'(k,t)
  =
  \e^{-\im k l} S(k) \bs{b}'(k,t)
  \\[2pt]
  =
  \frac{2 \e^{-\im k l} }{(\omega_1^{} - \omega_2^{}) (\e^{\im k l} - \e^{-\im k l})} \bigg[
  \begin{pmatrix}
    1
    \\
    \omega_2^{} / \im k
  \end{pmatrix}
  (\omega_1^{} - \theta k^2) \tilde{g}^{(2)}(t)
  -
  \begin{pmatrix}
    1
    \\
    \omega_1^{} / \im k
  \end{pmatrix}
  (\omega_2^{} - \theta k^2) \tilde{g}^{(1)}(t)
  \bigg]
\end{multline*}
(see Eq.~\eqref{eq:SaSbPrimeCaseOneM}).

In order to obtain a series representation, the paths of integration $\pd D_0^\pm$ must be deformed back to the real axis. Note that the branch points $k = (\pm \sqrt{\alpha \tau} \pm \sqrt{\alpha \tau - \mu^2})/\mu^2$ are removable. By deforming the paths, avoiding poles located on the real axis, and taking into account opposite directions of integration, we obtain contours around the poles. As a result of integration along the contours, the representation~\eqref{eq:SolutionCaseTwoIntegral} is transformed into the series
\begin{equation*}
  \bs{u}(x,t)
  =
  \frac{\pi \im}{2 \pi} \sum_{n \in \Z}^{} \Res_{k = k_n^{}}^{} \left[ \e^{\im k x} S(k) \bs{a}'(k,t) \right]
  + \frac{\pi \im}{2 \pi} \sum_{n \in \Z}^{} \Res_{k = k_n^{}}^{} \left[ \e^{\im k x} \e^{-\im k l} S(k) \bs{b}'(k,t) \right],
\end{equation*}
where $k_n^{} = \pi n / l$ are the poles.

Using the equalities $\e^{\pm \im k_n^{} l} = (-1)^n$, $[(\e^{\im k l} - \e^{-\im k l})_k']|_{k = k_n^{}}^{} = 2 \im l (-1)^n$ and $\omega_1^{}(k_n^{}) \omega_2^{}(k_n^{}) = \beta k_n^2$, we see that the residues are equal to
\begin{multline*}
  \Res_{k = k_n^{}}^{} \left[ \e^{\im k x} S(k) \bs{a}'(k,t) \right]
  =
  \Res_{k = k_n^{}}^{} \left[ \e^{\im k x} \e^{-\im k l} S(k) \bs{b}'(k,t) \right]
  \\[2pt]
  =
  \frac{\e^{\im k_n^{} x}}{\im l (\omega_{1,n}^{} - \omega_{2,n}^{})}
  \begin{pmatrix}
    (\omega_{1,n}^{} - \theta k_n^2) \tilde{g}^{(2)}(t) - (\omega_{2,n}^{} - \theta k_n^2) \tilde{g}^{(1)}(t)
    \\[5pt]
    - \im k_n^{} \big[ (\beta  -\theta \omega_{2,n}^{}) \tilde{g}^{(2)}(t) - (\beta  -\theta \omega_{1,n}^{}) \tilde{g}^{(1)}(t) \big]
  \end{pmatrix},
  \quad
  n \neq 0,
\end{multline*}
and
\begin{equation*}
  \Res_{k = 0}^{} \left[ \e^{\im k x} S(k) \bs{a}'(k,t) \right]
  =
  \Res_{k = 0}^{} \left[ \e^{\im k x} \e^{-\im k l} S(k) \bs{b}'(k,t) \right]
  =
  \frac{1}{\im l}
  \begin{pmatrix}
    \tilde{g}
    \\[5pt]
    0
  \end{pmatrix},
\end{equation*}
where $\omega_{m,n}^{} \equiv \omega_m^{}(k_n^{})$, and
\begin{equation*}
  \tilde{g}(t)
  =
  \int_0^tg(s) \diff s.
\end{equation*}

As a result, we conclude that the solution is given by
\begin{multline}
  \label{eq:SolutionCaseTwoSeries}
  \bs{u}(x,t)
  \equiv
  \begin{pmatrix}
    e(x,t) \\[2pt] q(x,t)
  \end{pmatrix}
  =
  \frac{1}{l}
  \begin{pmatrix}
    \tilde{g}(t) \\[3pt] 0
  \end{pmatrix}
  \\[2pt]
  + \frac{2}{l} \sum_{n=1}^\infty \frac{1}{\omega_{1,n}^{} - \omega_{2,n}^{}}
  \begin{pmatrix}
    \big[ (\omega_{1,n}^{} - \theta k_n^2) \tilde{g}^{(2)}(t) - (\omega_{2,n}^{} - \theta k_n^2) \tilde{g}^{(1)}(t) \big] \cos(k_n^{} x)
    \\[3pt]
    k_n^{} \left[ (\beta  -\theta \omega_{2,n}^{}) \tilde{g}^{(2)}(t) - (\beta  -\theta \omega_{1,n}^{}) \tilde{g}^{(1)}(t) \right] \sin(k_n^{} x)
  \end{pmatrix},
\end{multline}
where $\beta = \alpha/\tau$, $\theta = \mu^2/\tau$, $k_n^{} = \pi n / l$, $\omega_{m,n}^{} \equiv \omega_m^{}(k_n^{})$ (see Eq.~\eqref{eq:omegaOneTwo}), $\tilde{g}^{(1,2)}(t)$ are given by Eq.~\eqref{eq:tildeg}.
This solution can also be easily obtained by the Fourier method with preliminary elimination of inhomogeneities in boundary conditions. In Ref.~\cite{Kovacs:2018IJHMT} the solution was obtained for a particular heat flux $g(t)$ of the form~\eqref{eq:g}.

\subsection{The initial state is non-zero and the boundaries are thermally insulating}

Here we assume that heat sources are absent ($f = 0$), the initial state is non-zero, and the boundaries are thermally insulating ($\gamma_0^{} = \gamma_l^{} = 0$ and $g = 0$, $h = 0$). In this case the equations are homogeneous:
\begin{subequations}
  \label{eq:IBVPEqsCaseThree}
  \begin{align}
    &\pd_t^{} e + \pd_x^{} q
      = 0,
    \quad
    x \in (0,l),
    \quad
    t > 0,
    \\
    &\tau \pd_t^{} q
      + q
      - \mu^2 \pd_x^2 q
      + \alpha \pd_x^{} e
      = 0,
    \quad
    x \in (0,l),
    \quad
    t > 0.
  \end{align}
\end{subequations}
the initial conditions are
\begin{equation}
  \label{eq:IBVPInitCondCaseThree}
  \left. e \right|_{t=0}^{}
  =
  \varphi(x),
  \quad
  \left. q \right|_{t=0}^{}
  =
  \psi(x),
\end{equation}
and the boundary conditions are
\begin{equation}
  \label{eq:IBVPBoundCondCaseThree}
  \left. q \right|_{x=0}^{}
  = 0,
  \quad
  \left. q \right|_{x=l}^{}
  = 0.
\end{equation}

The solution to the problem~\eqref{eq:IBVPEqsCaseThree}--\eqref{eq:IBVPBoundCondCaseThree} is given by
\begin{multline}
  \label{eq:SolutionCaseThreeIntegral}
  \bs{u}(x,t)
  =
  \frac{1}{2\pi} \int_{-\infty}^\infty \e^{\im k x} \e^{-\varLambda(k) t} \hat{\bs{u}}_0^{}(k) \diff k
  \\
  - \frac{1}{2\pi} \int_{\pd D_0^{+}}^{} \e^{\im k x} S(k) \bs{a}'(k,t) \diff k
  - \frac{1}{2\pi} \int_{\pd D_0^{-}}^{} \e^{\im k x} \e^{-\im k l} S(k) \bs{b}'(k,t) \diff k.
\end{multline}
where
\begin{subequations}
  \label{eq:SaSbPrimeCaseThree}
  \begin{align}
    \label{eq:SaPrimeCaseThree}
    S(k) \bs{a}'(k,t)
    &=
      \frac{1}{\omega_1^{} - \omega_2^{}} \bigg[
      -
      \begin{pmatrix}
        1
        \\
        \omega_1^{} / \im k
      \end{pmatrix}
    \im k (\beta - \theta \omega_1^{}) \tilde{e}_0^{(1)}(0,t)
    +
    \begin{pmatrix}
      1
      \\
      \omega_2^{} / \im k
    \end{pmatrix}
    \im k (\beta - \theta \omega_2^{}) \tilde{e}_0^{(2)}(0,t)
    \bigg],
    \\[4pt]
    \label{eq:SbPrimeCaseThree}
    S(k) \bs{b}'(k,t)
    &=
      \frac{1}{\omega_1^{} - \omega_2^{}} \bigg[
      -
      \begin{pmatrix}
        1
        \\
        \omega_1^{} / \im k
      \end{pmatrix}
    \im k (\beta - \theta \omega_1^{}) \tilde{e}_0^{(1)}(l,t)
    +
    \begin{pmatrix}
      1
      \\
      \omega_2^{} / \im k
    \end{pmatrix}
    \im k (\beta - \theta \omega_2^{}) \tilde{e}_0^{(2)}(l,t)
    \bigg],
  \end{align}
\end{subequations}
and
\begin{subequations}
  \label{eq:eePrimeCaseThree}
  \begin{align}
    \label{eq:eOnePrimeCaseThree}
    -\im k (\beta - \theta \omega_1^{})
    \begin{pmatrix}
      \tilde{e}_0^{(1)}(0,t)
      \\[3pt]
      \tilde{e}_0^{(1)}(l,t)
    \end{pmatrix}
    =
    \frac{1}{\e^{\im k l} - \e^{-\im k l}} \Bigg[
    &
    \begin{pmatrix}
      \e^{\im k l}
      \\[3pt]
      1
    \end{pmatrix}
    V_1^{}(k,t)
    +
    \begin{pmatrix}
      \e^{-\im k l}
      \\[3pt]
      1
    \end{pmatrix}
    V_1^{}(-k,t)
    \Bigg],
    \\[4pt]
    \label{eq:eTwoPrimeCaseThree}
    \im k (\beta - \theta \omega_2^{})
    \begin{pmatrix}
      \tilde{e}_0^{(2)}(0,t)
      \\[3pt]
      \tilde{e}_0^{(2)}(l,t)
    \end{pmatrix}
    =
    \frac{1}{\e^{\im k l} - \e^{-\im k l}} \Bigg[
    &
    \begin{pmatrix}
      \e^{\im k l}
      \\[3pt]
      1
    \end{pmatrix}
    V_2^{}(k,t)
    +
    \begin{pmatrix}
      \e^{-\im k l}
      \\[3pt]
      1
    \end{pmatrix}
    V_2^{}(-k,t)
    \Bigg].
  \end{align}
\end{subequations}

Deforming the paths $\pd D_0^\pm$ back to the real axis, we see that the representation~\eqref{eq:SolutionCaseThreeIntegral} is transformed into the representation
\begin{multline}
  \label{eq:SolutionCaseThreeIntegralSeries}
  \bs{u}(x,t)
  =
  \frac{1}{2\pi} \int_{-\infty}^\infty \e^{\im k x} \e^{-\varLambda(k) t} \hat{\bs{u}}_0^{}(k) \diff k
  \\
  - \frac{1}{2\pi} \pvint_{-\infty}^\infty \e^{\im k x} S(k) \bs{a}'(k,t) \diff k
  + \frac{\pi \im}{2 \pi} \sum_{n \in \Z}^{} \Res_{k = k_n^{}}^{} \left[ \e^{\im k x} S(k) \bs{a}'(k,t) \right]
  \\
  + \frac{1}{2\pi} \pvint_{-\infty}^\infty \e^{\im k x} \e^{-\im k l} S(k) \bs{b}'(k,t) \diff k
  + \frac{\pi \im}{2 \pi} \sum_{n \in \Z}^{} \Res_{k = k_n^{}}^{} \left[ \e^{\im k x} \e^{-\im k l} S(k) \bs{b}'(k,t) \right],
\end{multline}
where $k_n^{} = \pi n / l$ are the poles.
It follows from Eqs.~\eqref{eq:SaSbPrimeCaseThree} and~\eqref{eq:eePrimeCaseThree} the relation
\begin{equation*}
  S(k) \bs{a}'(k,t)
  - \e^{-\im k l} S(k) \bs{b}'(k,t)
  =
  S(k) \bs\varPhi(k)
  =
  \e^{-\varLambda(k) t} \hat{\bs{u}}_0^{}(k).
\end{equation*}
Hence the difference of the two principal-value integrals in Eq.~\eqref{eq:SolutionCaseThreeIntegralSeries} is equal to the first integral, and, therefore, the net contribution of the three integrals to the solution is equal to zero. Thus the representation~\eqref{eq:SolutionCaseThreeIntegralSeries} takes the form
\begin{equation*}
  \bs{u}(x,t)
  =
  \frac{\im}{2} \sum_{n \in \Z}^{} \Res_{k = k_n^{}}^{} \left[ \e^{\im k x} S(k) \bs{a}'(k,t) \right]
  + \frac{\im}{2} \sum_{n \in \Z}^{} \Res_{k = k_n^{}}^{} \left[ \e^{\im k x} \e^{-\im k l} S(k) \bs{b}'(k,t) \right].
\end{equation*}
Using the equalities $\e^{\pm \im k_n^{} l} = (-1)^n$ and $[(\e^{\im k l} - \e^{-\im k l})_k']|_{k = k_n^{}}^{} = 2 \im l (-1)^n$ we see that the residues are equal to
\begin{multline*}
  \Res_{k = k_n^{}}^{} \left[ \e^{\im k x} S(k) \bs{a}'(k,t) \right]
  =
  \Res_{k = k_n^{}}^{} \left[ \e^{\im k x} \e^{-\im k l} S(k) \bs{b}'(k,t) \right]
  \\[2pt]
  =
  \frac{\e^{\im k_n^{} x}}{2 \im l} \left\{
    S(k_n^{}) \big[
      \bs\varPhi(k_n^{}) + \bs\varPhi(-k_n^{})
    \big]
  \right\}
  \\[2pt]
  =
  \frac{\e^{\im k_n^{} x}}{2 \im l (\omega_{1,n}^{} - \omega_{2,n}^{})}
  \begin{pmatrix}
    \varPhi_1^{}(k_n^{}) + \varPhi_1^{}(-k_n^{}) + \varPhi_2^{}(k_n^{}) + \varPhi_2^{}(-k_n^{})
    \\[2pt]
    \big\{ \omega_{1,n}^{} \left[ \varPhi_1^{}(k_n^{}) + \varPhi_1^{}(-k_n^{}) \right] + \omega_{2,n}^{} \left[ \varPhi_2^{}(k_n^{}) + \varPhi_2^{}(-k_n^{}) \right] \big\} / \im k_n^{}
  \end{pmatrix},
\end{multline*}
where $\omega_{m,n}^{} \equiv \omega_m^{}(k_n^{})$.
The functions $\varPhi_{1,2}^{}$ are given by Eq.~\eqref{eq:varPhi}, hence
\begin{multline*}
  \varPhi_1^{}(k_n^{}) + \varPhi_1^{}(-k_n^{})
  + \varPhi_2^{}(k_n^{}) + \varPhi_2^{}(-k_n^{})
  \\[2pt]
  \shoveleft{=
  (\omega_{1,n}^{} \e^{-\omega_{2,n}^{} t}  - \omega_{2,n}^{} \e^{-\omega_{1,n}^{} t}) \left[ \hat{\varphi}(k_n^{}) + \hat{\varphi}(-k_n^{}) \right]}
  \\[2pt]
  \shoveright{+ \im k_n^{} (\e^{-\omega_{1,n}^{} t} - \e^{-\omega_{2,n}^{} t}) \left[ \hat{\psi}(k_n^{}) - \hat{\psi}(-k_n^{}) \right]}
  \\[2pt]
  =
  2 \left[ (\omega_{1,n}^{} \e^{-\omega_{2,n}^{} t}  - \omega_{2,n}^{} \e^{-\omega_{1,n}^{} t}) \varphi_n^{}
    + k_n^{} (\e^{-\omega_{1,n}^{} t} - \e^{-\omega_{2,n}^{} t}) \psi_n^{} \right],
  \quad
  n \in \Z \backslash \{0\},
\end{multline*}
and
\begin{multline*}
  \frac{\omega_{1,n}^{} \left[ \varPhi_1^{}(k_n^{}) + \varPhi_1^{}(-k_n^{}) \right] + \omega_{2,n}^{} \left[ \varPhi_2^{}(k_n^{}) + \varPhi_2^{}(-k_n^{}) \right]}{\im k_n^{}}
  \\[2pt]
  \shoveleft{=
  \im k_n^{} \beta (\e^{-\omega_{1,n}^{} t} - \e^{-\omega_{2,n}^{} t}) \left[ \hat{\varphi}(k_n^{}) + \hat{\varphi}(-k_n^{}) \right]}
  \\[2pt]
  \shoveright{+ (\omega_{1,n}^{} \e^{-\omega_{1,n}^{} t} - \omega_{2,n}^{} \e^{-\omega_{2,n}^{} t}) \left[ \hat{\psi}(k_n^{}) - \hat{\psi}(-k_n^{}) \right]}
  \\[2pt]
  =
  2 \im \left[ k_n^{} \beta (\e^{-\omega_{1,n}^{} t} - \e^{-\omega_{2,n}^{} t}) \varphi_n^{}
  - (\omega_{1,n}^{} \e^{-\omega_{1,n}^{} t} - \omega_{2,n}^{} \e^{-\omega_{2,n}^{} t}) \psi_n^{} \right],
  \quad
  n \in \Z \backslash \{0\},
\end{multline*}
where
\begin{equation*}
  \begin{aligned}
    \label{eq:PhinPsinCaseThree}
    \hat{\varphi}(k_n^{}) + \hat{\varphi}(-k_n^{})
    =
    l \varphi_n^{},
    \quad
    \varphi_n^{}
    =
    \frac{2}{l} \int_0^l \varphi(x) \cos(k_n^{} x) \diff x,
    \quad
    \varphi_0^{}
    =
    \frac{1}{l} \int_0^l \varphi(x) \diff x,
    \\[2pt]
    \hat{\psi}(k_n^{}) - \hat{\psi}(-k_n^{})
    =
    - \im l \psi_n^{},
    \quad
    \psi_n^{}
    =
    \frac{2}{l} \int_0^l \psi(x) \sin(k_n^{} x) \diff x.
  \end{aligned}
\end{equation*}

As a result, we conclude that the solution is given by
\begin{multline*}
  \bs{u}(x,t)
  \equiv
  \begin{pmatrix}
    e(x,t) \\[2pt] q(x,t)
  \end{pmatrix}
  =
  \begin{pmatrix}
    \varphi_0^{} \\[3pt] 0
  \end{pmatrix}
  + \sum_{n=1}^\infty \frac{1}{\omega_{1,n}^{} - \omega_{2,n}^{}}
  \times
  \\[4pt]
  \times
  \begin{pmatrix}
    \big[
      (\omega_{1,n}^{} \e^{-\omega_{2,n}^{} t}  - \omega_{2,n}^{} \e^{-\omega_{1,n}^{} t}) \varphi_n^{}
      + k_n^{} (\e^{-\omega_{1,n}^{} t} - \e^{-\omega_{2,n}^{} t}) \psi_n^{}
    \big] \cos(k_n^{} x)
    \\[4pt]
    \big[
      - k_n^{} \beta (\e^{-\omega_{1,n}^{} t} - \e^{-\omega_{2,n}^{} t}) \varphi_n^{}
      + (\omega_{1,n}^{} \e^{-\omega_{1,n}^{} t} - \omega_{2,n}^{} \e^{-\omega_{2,n}^{} t}) \psi_n^{}
    \big] \sin(k_n^{} x)
  \end{pmatrix},
\end{multline*}
where $k_n^{} = \pi n / l$, $\omega_{m,n}^{} \equiv \omega_m^{}(k_n^{})$ (see Eq.~\eqref{eq:omegaOneTwo}).
This solution can also be easily obtained by the Fourier method.

\renewcommand{\bibname}{References}

\end{document}